\documentstyle[amscd,amssymb,12pt]{amsart}
\input xypic \xyoption{curve}
\input epsf.sty

\def\hcorrection#1{\advance\hoffset by #1 }
\def\vcorrection#1{\advance\voffset by #1 }

\vcorrection{-.7in}
\hcorrection{-.5in}

\newcommand{\C}[1]{{\cal#1}} % Calligraphic
\newcommand{\D}[1]{{\Bbb#1}} % Blackboard (for Z, Q, R, C, etc.)
\theoremstyle{plain}
\newtheorem{th}{Theorem}[section]
\newtheorem{cor}{Corollary}[section]
\newtheorem{lem}{Lemma}[section]
\newtheorem{prop}{Proposition}[section]

\theoremstyle{definition}
\newtheorem{defin}{Definition}[section]

\theoremstyle{definition}

\theoremstyle{remark}
\newtheorem{rem}{Remark}[section]

\numberwithin{equation}{section}

\begin{document}

\pagestyle{plain}
\addtolength{\footskip}{.3in}

% *****************  Title, Abstract & Introduction *******************
% *********************************************************************
%                      		Title page
% *********************************************************************
%Topmatter
\title{The Nonabelian Bar Resolution$^1$}
\author{Lucian M. Ionescu}
\address{Department of Mathematical Sciences\\
	University of Northern Colorado\\
        Greeley, Colorado 80639}
\email{lmiones@@fisher.unco.edu}
%\thanks{...} 
\date{01/01; $^1$ presented at the AMS meeting - Lawrence, KS, March 30-31, 2001}
%End topmatter

%=======================================================================
%                      		Abstract
%=======================================================================
\begin{abstract}
We develop the theory of parity quasi-complexes (PQC),
preparing the set up for defining derived functors using resolutions
in the nonabelian case.

A homotopy structure on the category of PQC is defined,
yielding a 2-category structure.
The nonabelian homology functor factors through 
the corresponding homotopy category.

Following the relative homological algebra approach,
resolutions are defined as PQC having parity contracting homotopies
in a suitable category.

A canonical non-abelian PQC resolution for groups is defined.
\end{abstract}

\maketitle
\tableofcontents

% *************************  end Header  *********************************

%=================================================================
%                   Introduction
%=================================================================
\section{Introduction}\label{S:intro}
The aim of this paper is to study parity quasicomplexes 
from the homotopical point of view, and to prepare the set up
for constructing derived functors on nonadditive functors using
parity resolutions. 
These are defined in homotopical terms,
as weak equivalences in the category of parity quasicomplexes.
The homotopy structure is defined such that weak equivalences
are quasi-isomorphisms with respect to the nonabelian homology functor.

The natural context is that of relative homological algebra,
adapted to the nonadditive case by considering parity quasi-complexes (PQC, \cite{PQC}).
The relevance of PQCs in relation with nonabelian cohomology was
noticed in \cite{PQC} (see also \cite{Hist}).

More specifically, PQC resolutions are defined as PQC having
(parity) contracting homotopies in a suitable category (\ref{D:ch}).
A non-abelian bar ``resolution'' for groups is defined \ref{D:barres}.
It has a canonical parity contracting homotopy (theorem \ref{T:barres}).

The usual notion of complex ($d^2=0$),
does not have a suitable analog in the nonabelian case,
and more general sequences are considered.

The notions of exactness and resolution,
equivalent in the context of a resolvent pair (\cite{ML}, p.265)
to the existence of a contracting homotopy (\cite{ML}, p.265),
is adapted to PQCs.
A PQC is exact if it has a parity contracting homotopy in the suitable ``larger'' category.

% The usual homological algebra arguments of factorization based on ...
If homological algebra is based on factorization arguments due to the $d^2=0$ condition
and lifting properties corresponding to a certain projective class
(e.g. the standard class determined by projective objects and epimorphisms),
the noncommutative case is ``simpler'', and depends only on the class
of weak equivalences (quasi-isomorphisms) in the sense of closed model categories.
In such a context, a cofibrant chain complex with a weak equivalence to a trivial complex,
is exactly a projective resolution (\cite{ht1}, p.113).
%Note that contracting homotopies are usually provided by splittings
%in a larger category.

We will construct a nonabelian bar resolution for a group $G$,
as an object of the category of $G$-groups,
having a (parity) contracting homotopy at the level of groups.
Forgetting the group structure, group extensions split as pointed sets.
The contracting homotopy allows the comparison of $n$-fold extensions of groups
with the bar resolution,
providing the correspondence with the group cohomology classes,
as in the abelian case.

% ******************************************************************************
%	Resolutions and homotopies
% ******************************************************************************
\section{Resolutions and homotopies}\label{S:resandhom}
We will fix the notation recalling some known facts,
relevant from the point of view of homotopical algebra.

In this section $\C{A}$ will denote an additive and exact category (e.g. $R$-mod),
with the associated category of complexes $Ch(\C{A})$
and homology functor $H_0$.
%
% an exact category; complexes; resolutions <=> quasi-isomorphisms
% contracting homotopy => quasi-iso; quasi-iso + splitting => c.h.

If $G$ is an object in $\C{A}$ then $G_\bullet=K(G,0)$ will denote the trivial complex
concentrated in degree 0 (the chain analog of Eilenberg-MacLane spaces). 

Let $X_\bullet\overset{\epsilon}{\to} G$ be an {\em augmented sequence}
of morphisms $d_n:X_n\to X_{n-1}$ in $\C{A}$, with $\epsilon\circ d_1=0$.
Then $X_\bullet\overset{\epsilon_\bullet}{\to} G_\bullet$ will denote the
corresponding chain map, and $(\tilde{X}_\bullet, \tilde{d}_\bullet)$
the corresponding sequence with $\tilde{X}_{-1}=G, \tilde{d}_{-1}=\epsilon$.

The chain map $\epsilon_\bullet$ is a {\em quasi-isomorphism}
iff it induces an isomorphism $H_\bullet(\epsilon_\bullet)$ in homology.
The chain map $\epsilon_\bullet$ has a
{\em contracting homotopy} $(\tau_\bullet, s_\bullet)$
iff there exists a chain map $\tau_\bullet:G_\bullet\to X_\bullet$
which is a section of $\epsilon_\bullet$,
i.e. $\epsilon_\bullet \tau_\bullet=id_{G_\bullet}$,
together with a homotopy $s_\bullet:1\to \tau_\bullet \epsilon_\bullet$
(\cite{ML}, p.41).

We will extend $s_\bullet$ to $\tilde{s}_\bullet$,
with the convention $\tilde{s}_{-1 }=\tau$.
\begin{lem}\label{L:qiso}
Let $\epsilon_\bullet:X_\bullet\to G_\bullet$ be a chain map in $Ch(\C{A})$,
where $G\bullet=K(G,0)$.

1) The augmented sequence $(X_\bullet, \epsilon)$
is a resolution of $G$ iff $\epsilon_\bullet$ is a quasi-isomorphism.

2) If $X_\bullet\overset{\epsilon}{\to}G$ is a projective resolution of $G$,
then $\epsilon:X_0\to G$ splits iff
$\epsilon_\bullet$ has a contracting homotopy.
\end{lem}
We will adapt the mechanism of relative homological,
where an exact sequence or a resolution (e.g. $R-mod$), has a contracting homotopy
in a ``larger category`` (e.g. $Ab$).
In the nonabelian case this property will replace the usual definition
of an exact sequence (see \ref{L:resol}).
As an example, in the abelian case,
the existence of a contracting homotopy in the ``base'' category ($\D{Z}-mod$),
allows the lifting of a given morphism to a chain transformation,
yielding the ``Comparison Theorem'' (\cite{ML}, T 6.1, p.87)
and the mechanism of derived functors based on resolutions.

% ****************************************************************
%		The category of parity quasicomplexes
% ****************************************************************

\section{The category of parity quasicomplexes}\label{S:catofpqc}
Let $\C{S}=Sets_*$ be the category of pointed sets, and
$\C{G}\C{S}$ the category of group objects in $\C{S}$,
i.e. the category of groups.
Denote by $\C{U}_\C{S}:\C{G}\C{S}\to \C{S}$ the forgetful functor and
$\C{F}_\C{S}:\C{S}\to \C{G}\C{S}$ the free object functor.
Then $(\C{U},\C{F})$ is a {\em resolvent pair of categories} \cite{ML}, p.265.

Let $G\in\C{G}$ be a fixed group and $\C{C}=G-\C{G}$ the category of $G$-groups,
i.e. the category with objects $(N, L)$, where $N$ is a group on which $G$ acts
through $L:G\to Aut(N)$, and with morphisms the $G$-equivariant group homomorphisms.
If $G$ is interpreted as a one-object category in the obvious way (\cite{Cat}),
then $\C{C}$ is the category of nonabelian modules in the categorical sense
($Hom(\C{C}_G, \C{G})$).
Let $U:\C{C}\to \C{S}$ the corresponding forgetful functor, and
$\C{F}:\C{S}\to \C{C}$ the free object functor.
\begin{lem}\label{L:resol}
A sequence of morphisms $X=(X_\bullet, d_\bullet)$ in $\C{C}$ is exact iff
$\C{U}(X)$ has a contracting homotopy.
\end{lem}
\begin{pf}
Note that $X$ is exact iff $U(X)$ is exact, and
that any exact sequence in $Sets_*$ splits.
Apply now the analog of 2 from lemma \ref{L:qiso}.
\end{pf}
\begin{defin}\label{D:PQC}
Let $Ch^\pm(\C{C})$ be the category of sequences of pairs of morphisms in $\C{C}$,
called {\em parity quasicomplexes} (PQC) \cite{PQC}:
\begin{gather}\xymatrix @C=3pc @R=3pc {
\dots
C_n \ar@/_/[r]_{\partial^-_n} \ar@/^/[r]^{\partial^+_n} \ar@{}[r]|{\partial_n \Uparrow}
& C_{n-1}
\dots 
} \label{D:pqc}\end{gather}
A {\em morphism of PQC} is a map $f\in Hom_0(X,Y)$
commuting with the parity differentials $\partial^\pm$.
\end{defin}
A PQC $(X_\bullet,\partial_\bullet)$ should rather be thought of
as sequence of ``2-morphisms'' $\partial_n:\partial^-_n\to \partial^+_n$:
%
% Define homology
%
\begin{defin}\label{D:homol}
The {\em $n^{th}$ homology space of a PQC},
is defined as the pointed set of cycles
$Z_n(X)$, the equalizer of $\partial_n^+$ and $\partial^-_n$,
modulo the equivalence relation defined as follows.
If $x, y \in X_n$, then $x\sim y$ iff there exists $c\in X_{n+1}$
such that $\partial^+c+x=y+\partial^-c$:
%
%			Hom(X,Y)
$$\diagram
H^{N}_n(X)=equal(\partial^+,\partial^-)/\sim & &
\dto_{x} \rto^{\partial^-c} & \dto^{y} \\
x \overset{c}{\sim} y & & \rto_{\partial^+c} &
\enddiagram$$
\end{defin}
Since a morphism of PQCs commutes with both quasidifferentials,
it induces as usual a morphism on homology.
\begin{prop}\label{P:homol1}
A morphism of PQCs, $f$, induces functorially a morphism of pointed sets
on the cohomology spaces, $H^{N}(f)$.
\end{prop}
In the above relation and in what follows,
the additive notation is used but without assuming commutativity.
Special care should be given when considering the opposite: $-(a-b)=b-a$.
To keep track of the order of the terms involved in a sum, categorical diagrams
are sketched to depict the required ``orientation''.

%
% 		Contracting homotopy is not enough
%
Parity quasicomplexes exhibit ``curvature'',
i.e. in general $\partial^2\ne 0$.
The usual notion of contracting homotopy is not sufficient,
as the following lemma demonstrates.
\begin{lem}\label{L:d2zero}
An augmented sequence of groups $(X_\bullet,\partial)$,
having a contracting homotopy $s$ in $Sets_*$
with image generating the groups $X_i$,
is a complex: $\partial^2=0$.
\end{lem}
\begin{pf}
The usual proof \cite{ML}, p.268, yields in the nonabelian case:
$$\partial\partial s+\partial=\partial+s \partial\partial.$$
Since it is a positive sequence, inductively follows $\partial^2s=0$,
and thus $\partial^2=0$.
\end{pf}
The notion of homotopy appropriate for PQC should not "mix"
the ``positive and negative boundaries'',
it should be a $\D{Z}_2$-graded version.

% *********************************************************************
%		PQC parity contracting Homotopies
%		Definitions: PCH, exact PQC/ resolution
% *********************************************************************
%				page 4
\section{The homotopy structure}
To define derived functors using resolutions, one needs that
homotopic maps induce the same map in homology (localization).
A useful homotopy should behave well with respect to
composition of morphisms of PQCs and functors (2-category structure).

To achieve this, it is natural to extend the principle of ``separation of signs''
to homotopies, and define a $\D{Z}_2$-graded analog.
\begin{defin}\label{D:ph}
Let $f, g:(X,\partial^\pm)\to (Y,\partial^\pm)$ be morphisms of PQCs.
$f$ is {\em parity homotopic} to $g$,
and denoted $f\overset{s}{\to}g$ iff for all $n$, $X_n$ is generated by the 
elements satisfying the following relations:
\begin{align}\label{E:ph}
f+s\partial^-&=\partial^+s+g\\
s\partial^+&=\partial^-s.
\end{align}
\end{defin}
\begin{rem}
The above addition is the ``natural addition'' of morphisms
($+_Y\circ (f\otimes g) \circ \Delta_X$).
The ``free addition'' (\cite{F2}, p.211), consisting in extending uniquely
a map defined on free generators (``through adjunction''), is used to
define the parity differentials, which satisfy the parity homotopy
relations only on generators.
This explains the courser relation defined above,
which is typical in the noncommutative case 
(see d.g. near-rings and pseudohomomorphisms: \cite{F1,F2}, \cite{L} p.313).

Note also that the homotopy structure defined by the requirement that the above 
relations hold globally, 
still enjoys the properties stated below (transitive, compatible with $H^N$, etc.).
\end{rem}
The homology functor (definition \ref{D:homol}) will factor through the
homotopy category (\cite{H}, p.126).
\begin{lem}\label{L:homol1}
Parity homotopic morphisms $f\to g:X\to Y$ induce the same map in homology:
$H^{N}(f)=H^{N}(g)$.
\end{lem}
\begin{pf}
If $x$ is a cycle, i.e. $\partial^-x=\partial^+x$, and $f\to g$, then:
$$\partial^-sx=s\partial^+x=s\partial^-x,$$
so that:
$$f(x)+\partial^-sx=\partial^+sx+g(x)$$
i.e. $f(x)$ and $g(x)$ are homologous cycles in $Y$.
\end{pf}
There is a natural ``vertical'' composition of homotopies.
\begin{lem}
The parity homotopy is transitive.
\end{lem}
\begin{pf}
If
$$f+s\partial^-=\partial^+s+g, \quad s\partial^+=\partial^-s$$
and
$$g+t\partial^-=\partial^+t+h, \quad t\partial^+=\partial^-t$$
then
$$f+(s+t)\partial^-=\partial^+s+g+t\partial^-=\partial^+(s+t)+h.$$
and
$$(s+t)\partial^+=\partial^-(s+t).$$
\end{pf}
Parity homotopy is distributive with respect to composition of morphisms.
\begin{lem}
If $f,g:X\to Y$ are parity homotopic $s:f\to g$, and $h:Y\to Z$ ($h:Z\to X$), 
then $hs:hf\to hg$ ($sh:fh\to gh$) is a parity homotopy between $hf$ and $hg$ 
($fh$ and $gh$).
\end{lem}
\begin{pf}
We will prove left distributivity. 
Left multiply equations \ref{E:ph} by $h$. Since $h$ commutes with the
parity quasidifferentials, the statement follows.
\end{pf}
The compatibility with composition of PQC morphisms follows.
\begin{cor}\label{L:scomp1}
Let $f,g:X\to Y$ and $f',g':Y\to Z$ be morphisms of PQCs.
If $f\to g$ and $f'\to g'$ then $f'f\to g'g$.
\end{cor}
\begin{pf}
The parity homotopy relation $\sim$ is left and right distributive
with respect to composition of PQC morphisms.
Then $f'f\sim g'f\sim g'g$.
\end{pf}
The above facts can be summarized as follows.
\begin{th}
$Ch^\pm(\C{C})$ is a 2-category, and the nonabelian homology functor
$H^{N}$ factors through the corresponding homotopy category.
\end{th}

%
% 		Homotopic functors - later!
%
%Let $\C{D}$ be an other category and $Ch^\pm(\C{D})$
%the associated category of PQCs.
%A functor $T:\C{C}\to \C{D}$ induces a functor $T^\#=Ch^\pm(T)$
%on the corresponding categories of PQCs.
%It is called a {\em homotopic functor} if it preserves the homotopy relation:
%if $f\sim g$ then $T^\#f\sim T^\#g$.

% f -> f# -> H(f)
% Homotopic functors induce ...

% Hom(.,N) and N\otimes. are homotopic functors?

%
\begin{defin}\label{D:ch}
A {\em parity contracting homotopy} (PCH) is a map $s\in End_1(X)$, such that $0\overset{s}{\to}id_X$:
\begin{alignat}{2}\label{E:pch}
&s\partial^- & =&\partial^+s+id_X\\
   &s\partial^+ & =&\partial^-s
\end{alignat}
A PQC $X$ is called {\em exact} if it has a parity contracting homotopy.
An augmented PQC $\epsilon:X_\bullet\to G$ is called a
{\em PQC-resolution of $G$} if it is an exact PQC.
\end{defin} % D:ch
%
%		lemma - s.c.h => trivial cohomology
%
\begin{cor}\label{C:ch}
An exact PQC has trivial cohomology.
\end{cor}
\begin{pf}
Using Lemma \ref{L:homol1}, yields $0=H^{N}(0_X)=H^{N}(id_X)=id_{H^{N}(X)}$.
\end{pf}
%		Remark
The direct shows that, in an exact PQC, every cycle is a boundary having 
a canonical decomposition into a ``positive'' and a ``negative '' part.
\begin{prop}
If $x$ is a cycle of an exact PQC, then it decomposes canonically
into a positive and a negative boundary:
$$x=\partial^+(-sx)-\partial^-(-s)x.$$
\end{prop}
\begin{pf}
As before, $\partial^-x=\partial^+x$ in combination with (\ref{E:pch}):
$$s\partial^-x=\partial^+sx+x, \quad s\partial^+x=\partial^-sx,$$
yields $\partial^-sx=\partial^+sx+x$. 
Rearranging the terms gives the above relation.
\end{pf}
We will prove in section \ref{S:barres} that the noncommutative bar resolution
(\S \ref{S:barres}) has a parity contracting homotopy.
%

% ****************************************************************
%			NA bar resolution
% ****************************************************************
%				page 5

\section{The nonabelian bar resolution}\label{S:barres}
The results are formulated for the category of groups
as a typical nonadditive category.
We are adapting the classical theory of derived functors based on
projective resolutions, following \cite{ML}.
We are using the notation from \cite{ML}, p.114,
with additive notation used for clarity.
The relative homological algebra facts are adapted 
to the noncommutative case.

The major difference is that the noncommutative case requires a
$\D{Z}_2$-grading of the usual bar resolution,
i.e. a ``separation'' of the positive and negative terms
(inverses, in the multiplicative case). 

The free constructions specified bellow refer to the category of
$G-groups$ \cite{F1}, 
consisting of groups together with a group action, and equivariant
group homomorphisms.
\begin{defin}\label{D:barres}
The {\em (nonabelian) bar resolution} of $G$ is the PQC
$(B_\bullet(G),\partial^\pm_\bullet)$:
\begin{gather}\xymatrix @C=3pc @R=3pc {
\dots
& B_2 \ar@/_/[r]_{\partial^-_2} \ar@/^/[r]^{\partial^+_2}
& B_1 \ar@{.>}[l]|{e_1} \ar@/_/[r]_{\partial^-_1} \ar@/^/[r]^{\partial^+_1}
& B_0 \ar@{.>}[l]|{e_0} \ar@/_/[r]_{0} \ar@/^/[r]^{\epsilon} 
& \ar@{.>}[l]|{e_{-1}} \D{Z} \to 0
}\label{D:ncbar}\end{gather}
% Example of a 2-cell
% \ar@/^1pc/[rr]^{(s,f)}  \ar@{}[rr]|{\Downarrow \gamma} \ar@/_1pc/[rr]_{(s',f')}
%
of $G$-groups $B_n(G)$ and $G$-equivariant group homomorphisms $\partial^\pm_n$,
together with the augmentation $\epsilon$, 
and the group homomorphisms $e_n$, defined below.
\end{defin}
$B_n=\C{F}(\C{U}(G^n))$ is the free $G$-group with ($G$-)generators $[x_1|...|x_n]$
all $n$-tuples of elements $x_1,...,x_n$ of $G$.

Operation on a generator with an element $x\in G$ yields an element $x[x_1|...|x_n]$ in $B_n$,
so $B_n$ may be described as the free group generated by all $x[x_1|...|x_n]$.
In particular, $B_0$ is the free $G$-group on one generator $[\ ]$,
so it is isomorphic to $\D{Z}G$, $B_{-1}=\D{Z}$ and $B_n=1$ for $n<-1$.

The group homomorphisms $e_n$ are defined by:
$$e_{-1}(1)=-[\ ]\ , \quad e_n(x[x_1|...|x_n])=-[x|x_1|...|x_n], \quad n\ge 0.$$
Similar to the abelian case (\cite{ML}, theorem 6.3, p.268),
the PQC structure maps are determined, if we require $e_n$ 
to be a parity contracting homotopy.
The negative sign included in the definition of $e_n$ will reverse the order of
odd simplicial maps in the structure formula of $\partial^-$.
%
% 	The Contracting Homotopy Theorem
%
\begin{th}\label{T:barres}
There are unique PQC structure morphisms $\partial^\pm$,
such that $e_\bullet$ is a parity contracting homotopy in $\C{G}$ for
$(B_\bullet,\partial_\bullet)$.
They satisfy the following structural equations on group generators:
\begin{align}
\partial^+_n=\sum_{\overset{0\le i\le n}{i\ even}}^{\to} \partial^i_n\qquad
\partial^-_n=\sum_{\overset{0\le i\le n}{i\ odd}}^{\leftarrow} \partial^i_n
\end{align}
\end{th}
\begin{pf}
Since the image of $e_i$ generates $B_{i+1}$,
the relations \ref{E:pch} inductively define $\partial^\pm$.

To prove the stated formulas, note that on group generators,
applying $e$ after a simplicial map $\partial^i$ yields $\partial^{i+1}s$:
$$\partial^{k+1}_{n+1} e_n=e_{n-1}\partial^k_n, k\ge 0,\qquad 
\partial^0_{n+1} e_n=-id_n.$$
Then, on $G$-group generators, we have:
$$\partial^+_{n+1}([x,y_1,...y_n])=-\partial^+_{n+1}e_n(x[y])=-[e_{n-1}\partial^-_n-id_n (x[y])],$$
and
\begin{alignat}{1}
e_{n-1}\partial^-_n(x[y])&=...+e\partial^3_n(x[y])+e\partial^1_n(x[y])\\
&=...+\partial^4_n e(x[y])+\partial^2_n e(x[y])\\
&=-(\partial^2_n+\partial^4_n+...)[x,y].
\end{alignat}
Since $id_n=\partial^0_n e_{n-1}$, the formula for $\partial^+_{n+1}$ is established.

The computations for $\partial^-_{n+1}$ are similar, and will be omitted.
\end{pf}
The above $G$-equivariant group homomorphisms $\partial^\pm$ (compare \cite{PQC}, p.6),
are expressed in terms of the usual simplicial maps defined on $G$-generators by:
\begin{gather}
\partial^0_n [x_1|...|x_n]=x_1[x_2|...|x_n],\qquad
\partial^n_n [x_1|...|x_n]=[x_1|...|x_{n-1}]\\
\partial^i_n [x_1|...|x_n]=[x_1|...|x_i x_{i+1}|...|x_n], \ i=1,...,n-1\notag
\end{gather}
In particular $s_0(x[\ ])=[x]$ and:
$$ 
\partial^+_1[x]=x[\ ], \quad \partial^-_1[x]=[\ ], \qquad
\partial^+_2[x|y]=x[y]+[x], \quad \partial^-_2[x|y]=[xy]
$$
\begin{rem}
We think of the sequence of pairs of maps as a ``horizontal'' (1-)sequence of
2-morphisms $\partial_n:\partial^-_n\to\partial^+_n$, in relation with the
globular approach to weak n-categories.
\end{rem}
\begin{defin}
The $G$-equivariant group homomorphism defined on generators by:
$$\partial_n([x_1,...,x_n])=\partial^+_n([x_1,...,x_n])-\partial^-_n([x_1,...,x_n]),$$
will be called the quasi-differential of the PQC $(X, \partial^\pm)$.
\end{defin}
%

% **************************************************************
%		PCHreier theory revisited
% **************************************************************
\section{Schreier theory revisited}\label{S:Schreier}
We reinterpret Schreier's theory of group extensions (\cite{CC}, p.198),
in terms of the underlying contracting homotopy.
The notation used is essentially that of \cite{ML} (4.3, 4.4, p.111).
\begin{gather}\diagram
\C{E} & 0 \rto^{} &
N \ar@{=}[d]_{id} \rto^{\chi} &
E \ar@{=}[d]_{id} \ar@{.>}[dl]_{v} \rto^{\sigma} &
G \ar@{=}[d]_{id} \ar@{.>}[dl]_{u} \rto^{} & 1 & \\
\C{E} & 0 \rto^{} &
N \rto^{\chi} &
E \rto^{\sigma} &
G \rto^{} & 1 & \\
\enddiagram\label{D:gext}\end{gather}
The elements of $B$ can be represented in a unique way as
$b=\chi(a)+u(x)$, where $x=\sigma(b)$.
Define $v(b)=a$, to have:
\begin{equation}\label{E:vext}
id_B=\chi\circ v+u\circ \sigma.
\end{equation}
Relabeling: $s_0=u, s_1=v, d_1=\sigma, d_2=\chi$, 
$s_\bullet$ is a contracting homotopy in $Sets_*$,
for the group extension $E_\bullet$: 
$$id_E\sim 0:\qquad d_0 s_0=id_G, \quad d_2 s_1+s_0 d_1=id_E, \quad s_1 d_2=id_N$$
\begin{rem}\label{R:adjext}
The contracting homotopy $s_\bullet$ defines a
``adjoint'' extension in $Sets_*$:
\begin{gather}\diagram
\C{E}^* &
0 \rto &
G                          \ar@<2pt>[r]^{u} &
E \ar@{.>}@<2pt>[l]^{\sigma} \ar@<2pt>[r]^{v} &
N \ar@{.>}@<2pt>[l]^{\chi} \rto & 1 & 
\enddiagram\label{D:adjext}\end{gather}
with a contracting homotopy $s^*_0=\chi, s^*_1=\sigma$
defined by the original extension.
\end{rem}
The meaning of the factor set $f(x,y)$ of the extension
is apparent when comparing the group extension with the bar resolution.

The correspondence between group extensions and
cohomology classes of the bar resolution is established
through $R-mod$ extensions (\cite{ML} T 6.2, p.121).

The general case requires defining the characteristic extension at 
the level of $R-groups$, where $R=Z(G)$ is the distributively
generated near-ring associated to $G$ (cite{F}),
and will be addressed elsewhere.

Consider the following diagram:
\begin{gather} \xymatrix @C=3.5pc @R=3.5pc {
B(G) &
B_3 \dto^{} \ar@/^/[r]^{\partial^+_3} \ar@/_/[r]_{\partial^-_3} &
B_2 \dto_{g_2} 
      \ar@/^/[r]^{\partial^+_2} \ar@/_/[r]_{\partial^-_2} &
B_1 \ar@{.>}[dl]_{\tilde{s}_1} \dto_{g_1}
      \rto^{ev} &
G \ar@{.>}[dl]_{\tilde{s}_0} \ar@{=}[d] \rto^{} & 0\\
\C{E} & 0 \rto^{} &
N \ar@<2pt>[r]^{d_2} & \ar@{.>}@<2pt>[l]^{s_1}
E \ar@<2pt>[r]^{d_1} & \ar@{.>}@<2pt>[l]^{s_0} 
G \rto^{} & 0
}\label{D:modext2}\end{gather}
where $\C{E}$ is the above group extension, and the top row
is a truncation of the bar resolution, 
``augmented'' with $ev$, obtained lifting the group
multiplication of $G$ to $B_1$, the free $G-group$ generated by $G$.
Left action of $G$ is inner conjugation.

There is a unique $G$-equivariant lift of the section $s$, denoted $\tilde{s}$
(see theorem below).
As in the classical case $g_2 e_0=s_0 g_1 (\partial^+_2-\partial^-_2)$
defines a $G$-equivariant chain map $d_2 g_2=g_1 \partial_2$, 
where $\partial_2$ is the $G$-equivariant morphism defined on generators
as $\partial^+_2-\partial^-_2$.
Then, on group generators:
$$g_2([x,y])=s(x)+s(y)-s(xy)=f(x,y)$$
is the factor set corresponding to the section $s$.
$g_2$ is a 2-cocycle relative to the
pseudoaction $L_x(e)=s(x)+e-s(x)$, induced by $G$ on $N$:
$$g_2\circ \partial^+_3=g_2\circ \partial^-_3,$$
relation which, in terms of the factor set $f$, 
takes the following form:
$$L_x(f(y,z))+f(x,yz)=f(x,y)+f(xy,z), x,y,z\in G.$$
It is well known that $L$ and $f$ are compatible in the following sense:
\begin{equation}\label{E:MC}
L_x L_y(n)+f(x,y)=f(x,y)+L_{xy}(n), \ n\in N,\ x,y\in G.
\end{equation}
The interpretation of this relation is given in the theorem below, 
which can be proved through a direct computation.
\begin{th}\label{T:MC}
With the above notations, the following are equivalent:

1) Equation \ref{E:MC} holds;
2) $g_1$ is $G$-equivariant;
3) $g_1\circ\partial_3\circ\partial_2=0$.
\end{th}

%=======================================================================
%                      Bibliography
%=======================================================================

% **********************************************************************

\end{document}